\documentclass[10pt,reqno]{amsart}
\usepackage{amssymb,mathrsfs,color}
\usepackage{pinlabel}

\usepackage{amsmath, amsfonts, amsthm, verbatim, amssymb}
\usepackage{epstopdf, mathtools}
\usepackage{color}
\usepackage{esint}
\usepackage{stmaryrd}
\usepackage{graphicx}
\usepackage{xcolor} 
\usepackage{tensor}
\usepackage{slashed}
\usepackage{cite}
\mathtoolsset{showonlyrefs}

\newcommand{\bs}[1]{\boldsymbol{#1}}

\newcommand{\cl}[1]{\mathcal{#1}}

\newcommand{\al}{\alpha}

\newcommand{\la}{\lambda}

\newcommand{\Om}{\Omega}

\newcommand{\p}{\partial}
\newcommand{\grad}{\nabla}

\usepackage[multiple]{footmisc}
\numberwithin{equation}{section}

\theoremstyle{remark}

\newcommand{\tr}{\mathrm{tr}\,}

\renewcommand{\div}{\mathrm{div}\,}
\newcommand{\rar}{\rightarrow}

\newcommand{\bbA}{\mathbb A}
\newcommand{\bbB}{\mathbb B}

\newcommand{\bbD}{\mathbb D}
\newcommand{\bbE}{\mathbb E}

\newcommand{\bbG}{\mathbb G}

\newcommand{\bbR}{\mathbb R}

\newcommand{\calL}{\mathcal L}

\newcommand{\calS}{\mathcal S}

\newcommand{\tens}{\otimes}

\newcommand{\Div}{\mbox{Div}\,}

\begin{document}

\title[Second-Gradient Incompressible Viscous Fluids]{Second-gradient models for incompressible viscous fluids and associated cylindrical flows} 
\author{C. Balitactac and C. Rodriguez}

\begin{abstract}
We introduce second-gradient models for incompressible viscous fluids, building on the framework introduced by Fried and Gurtin.  We propose a new and simple constitutive relation for the hyperpressure to ensure that the models are both physically meaningful and mathematically well-posed. The framework is further extended to incorporate pressure-dependent viscosities. We show that for the pressure-dependent viscosity model, the inclusion of second-gradient effects guarantees the ellipticity of the governing pressure equation, in contrast to previous models rooted in classical continuum mechanics. The constant viscosity model is applied to steady cylindrical flows, where explicit solutions are derived under both strong and weak adherence boundary conditions. In each case, we establish convergence of the velocity profiles to the classical Navier-Stokes solutions as the model's characteristic length scales tend to zero.

\end{abstract}

\maketitle

\section{Introduction}

\subsection{Higher-gradient continua} This work introduces \textit{second-gradient} models for incompressible viscous fluids including those with pressure-dependent viscosities. Here, ``second-gradient" refers to the highest order of spatial derivatives of the velocity appearing in the internal power expenditures. The general model introduced in \cite{FriedGurtin06} contains an internal ambiguity (outlined in Section 1.2) that restricts its predictive capability. The aim of this paper is to resolve this issue and to develop physically meaningful and mathematically well-posed theories that extend the classical incompressible Navier–Stokes model into regimes where it breaks down, including \textit{small-length scale flows} and \textit{high-pressure flows} exhibiting pressure-dependent viscosities; see \cite{Koplik88, Koplik_1995, XiChwang03, OkamuraHeyes04, Travis00, Travis97} and Bridgman's Nobel Prize winning work \cite{Bridgman1926, Bridgman1931, BridgmanNobelPrize}.

In \textit{classical continuum mechanics}, only the \textit{first gradient} of the spatial velocity appears in the internal power expended during the motion, and models do not possess internal length scales. The idea of continua with internal power or work expenditures depending on second- and higher-order spatial derivatives of the motion was first conceived by Piola in 1846 \cite{Piola1846Book, dellIsola15}. Later, the Cosserat brothers \cite{Cosserats1909} introduced additional degrees of freedom at each material point, represented by an orthonormal frame parameterizing a rigid microstructure. However, major advances in generalized continuum theories incorporating small-length scale effects did not occur until the latter half of the 20th century, through foundational work by Toupin \cite{Toupin62, Toupin64}, Green and Rivlin \cite{GreenRivlin64a, GreenRivlin64b}, Mindlin \cite{Mindlin64a, Mindlin1965}, Mindlin and Eshel \cite{MindlinEshel1968}, and Germain \cite{Germain73a, Germain73b}. These efforts culminated in general theories of gradient and micromorphic continua characterized by possibly multiple intrinsic length scales (see, e.g., \cite{Askes2011, maugin2011, maugin2013, dellIsola17, Maugin17Book, dellIsola2020higher, bertram2020}).

\subsection{Fried-Gurtin model of second-gradient incompressible linear viscous fluids} The models developed here builds upon a general framework for incompressible linear viscous fluids introduced in \cite{FriedGurtin06}. In particular, an additional third-order \textit{hyperstress} tensor $\bbG$ is introduced alongside the classical Cauchy stress tensor $\boldsymbol T$, with $\bbG$ dual to the second gradient of the spatial velocity $\boldsymbol v$ in the internal power expenditure\footnote{For a discussion of our notation, see Section 2.1.}
\begin{align}
	\quad \forall \cl P_t \subseteq \cl B_t, \quad 
	P_{\mathrm{int}}(\mathcal P_t, \boldsymbol v) = \int_{\mathcal P_t} (\boldsymbol T \cdot \nabla \boldsymbol v + \bbG \cdot \nabla \nabla \boldsymbol v)\, dv.
\end{align}
The resulting governing equations are
\begin{align}
	\boldsymbol 0 = \mathrm{div}\,(\boldsymbol T - \mathrm{div}\,\bbG) + \rho (\boldsymbol b + \boldsymbol b^{\mathrm{in}}), \quad \mathrm{div}\,\boldsymbol v = 0, \label{eq:I1}
\end{align}
where $\rho$ is the constant mass density, $\boldsymbol b$ {a} specific external body force, and $\boldsymbol b^{\mathrm{in}}$ the specific inertial body force; see Section 2.2. For incompressible linear viscous fluids with constant viscosity $\mu > 0$, the classical form $\boldsymbol T = -p \boldsymbol I + 2 \mu \boldsymbol D$ where $p$ is the mean normal stress (or pressure) is adopted. The hyperstress $\bbG$ takes the form
\begin{align}
	\bbG = -\frac{1}{2}[\boldsymbol I \otimes \boldsymbol \pi + (\boldsymbol I \otimes \boldsymbol \pi)^{T(2,3)}] + \bbG_0,
\end{align}
where $\boldsymbol \pi$ is an indeterminate vector field known as the \textit{hyperpressure}, and $\bbG_0$ is a linear isotropic function of $\nabla \nabla \boldsymbol v$ possessing the same symmetries and trace properties as $\nabla\nabla \boldsymbol{v}$, determined by three intrinsic length scales. This structure introduces additional degrees of modeling flexibility but also leaves a critical gap: the hyperpressure $\boldsymbol \pi$ cannot be determined directly from the field equations.

The absence of a constitutive relation for the hyperpressure limits the predictive power of the general second-gradient model for incompressible viscous fluids. The model requires additional assumptions to be physically complete since $\boldsymbol \pi$ appears in both the tractions and \textit{hypertractions} needed to specify boundary conditions; see \eqref{eq:4}, \eqref{eq:5} and Section 2.6. The objective of this work is to propose a specific form and physical interpretation of the hyperpressure (see \eqref{eq:piintro} below), and to formulate a corresponding model for second-gradient incompressible fluids with pressure-dependent viscosities. Filling this gap enables the development of physically meaningful and mathematically well-posed theories that extend classical viscous flow models to small-length scale and high-pressure regimes where the classical incompressible Navier–Stokes model fails.

Although the lack of specification for the hyperpressure presents a major theoretical limitation, it has not prevented the study of the general second-gradient framework in various settings. Many studies have employed either strong adherence boundary conditions, where the velocity and its normal derivative vanish at the boundary, or weak formulations based on the principle of virtual power, which only require knowledge of the ``effective pressure" $p - \mathrm{div}\,\boldsymbol \pi$ \cite{friedcontinuum2008, giusterithree-dimensional2010, giusterinonsimple2011, giusterislender-body2014, bodnarexistence2021, straughanthermal2023, giantesiothermal2024}. In these contexts, the need for an explicit constitutive relation for $\boldsymbol \pi$ is bypassed, allowing for either computing explicit solutions for concrete problems or establishing well-posedness results for the governing field equations without resolving the nature of the hyperpressure. Thus, while the general Fried-Gurtin framework has been widely employed, a physically grounded and predictive specification of the hyperpressure remains an open question, one that this paper seeks to address.  

\subsection{Outline} A brief outline of this work is as follows. Section 2 introduces the mathematical formulation of the second-gradient model for incompressible linear viscous fluids, including the governing equations, complementary boundary conditions, and a constitutive relation for the hyperpressure. In particular, motivated by recent arguments from \cite{krawietz2022surface}, we propose in Section 2.4 that 
\begin{align}
\bs \pi = \ell_1^2 \nabla p, \label{eq:piintro} 
\end{align}
where $\ell_1^2$ is a certain linear combination of the three independent length scales of the general second-gradient model. We also propose an additional boundary condition necessary to uniquely determine the pressure based on \eqref{eq:piintro}; see \eqref{eq:pressureboundary}. Section 3 extends the model to incorporate pressure-dependent viscosities and demonstrates how \eqref{eq:piintro} \textit{ensures ellipticity of the pressure equation}, addressing shortcomings of earlier models rooted in classical continuum mechanics. Section 4 applies the constant viscosity model to steady Poiseuille flow in a cylindrical tube, deriving explicit solutions under both strong and weak adherence boundary conditions and establishing convergence to the classical parabolic profile as the model's intrinsic length scales tend to zero. Section 5 examines rotating {Taylor-Couette} flow in a cylinder, derives explicit expressions for the velocity and pressure fields for strong and weak adherence boundary conditions, and shows convergence of the velocity profiles to the classical counterpart as the model's intrinsic length scales tend to zero.

\subsection*{Acknowledgments} The authors thank Jeremy Marzuola for his help generating numerical approximations of the pressure in the final section. C. B. and C. R. gratefully acknowledge support from NSF DMS-2307562, and C. B. also gratefully acknowledges support from NSF RTG DMS-2135998.

\section{A Second-gradient Incompressible Viscous Fluid Model}

\subsection{Preliminaries} In what follows, we use standard tensor notation, algebra, and calculus. In particular, bold lower-case font, $\bs a$, is used to denote vectors in $\bbR^3$, bold upper-case font, $\bs A$, is used to denote second-order tensors on $\bbR^3$, and bold blackboard font, $\mathbb A$, is used to denote third-order tensors on $\bbR^3$. The divergence of a second- or third-order tensor will always be taken with respect to the last index. For third-order tensors, $\mathbb A^{T(i,j)}$ and $\tr_{(i,j)}\bbA$ denotes taking the transposition and trace, respectively, with respect to the $i$th and $j$th indices. Finally, the gradient of a second-order tensor field $\bs A$ is the third-order tensor field, $\nabla \bs A$ given by
\begin{align}
	\nabla \bs A =  \Bigl (\frac{\p}{\p \theta^i}\bs A\Bigr ) \tens  \bs g^i,
\end{align}
where $\{\theta^i\}$ is a general curvilinear system of coordinates for three-dimensional Euclidean space $\bbE^3$ with dual basis $\{ \bs g^i \}$. 
We abuse notation slightly by denoting the scalar products between the various types of tensors in the same way: in Cartesian coordinates
\begin{align}
 \bs A \cdot \bs B := \sum_{i,j = 1}^3 A_{ij}B_{ij}, \quad \bbA \cdot \bbB := \sum_{i,j,k = 1}^3 A_{ijk} B_{ijk}.	
\end{align}
The corresponding Frobenius norms are then defined via 
\begin{align}
	|\bs A| := (\bs A \cdot \bs A)^{1/2}, \quad |\bbA| = (\bbA \cdot \bbA)^{1/2}. 
\end{align}

We work in the spatial (or Eulerian) description of the fluid. In what follows, $\mathcal B_t \subseteq \bbE^3$ is the current configuration of the body at time $t$, $\bs x \in \cl B_t$, and $\cl P_t \subseteq \cl B_t$ denotes a generic part of $\cl B_t$. We assume throughout that $\cl B_t$ and $\cl P_t$ have smooth boundaries. Finally, the material time derivative of a scalar or tensor valued variable $\bs f$ defined on $\cl B_t$ is denoted $\dot{\bs f}$. 

The \textit{mass density} is denoted by $\rho(\bs x,t)$, and the \textit{velocity} field is denoted by $\bs v(\bs x,t)$, with $\bs x \in \cl B_t$. The \textit{stretching} and \textit{spin} are denoted, respectively, by 
\begin{align}
	\bs D = \frac{1}{2}( \nabla \bs v + [\nabla \bs v]^T ), \quad \bs W = \frac{1}{2}(\nabla \bs v - [\nabla \bs v]^T ).
\end{align}

The local equation expressing conservation of mass is given by 
\begin{align}
	\dot \rho + \rho (\mathrm{div}\, \bs v) = 0, \quad \mbox{on } \cl B_t. \label{eq:massbalanceog}
\end{align}
In this work, we consider only incompressible fluids with constant density for which $\div \bs v = 0$ and \eqref{eq:massbalanceog} is automatically satisfied. 

\subsection{Fried-Gurtin model of non-simple materials}
 In \cite{FriedGurtin06}, Fried and Gurtin apply the principle of virtual power to each part of the body, first developed in \cite{Germain73a, Germain73b}, to derive a model for materials along with complementary boundary conditions that incorporate higher-order spatial derivatives into the equations of motion.

More precisely, it is assumed that at some fixed time $t$, $\cl B_t$, $\bs v$, $\dot{\bs v}$, $\rho$, and a non-inertial body force $\bs b$ are known. For a given vector field $\bs u$ on $\cl B_t$ (the \textit{test velocity} or \textit{virtual velocity}), the \textit{virtual power of the internal forces} for $\cl P_t \subseteq \cl B_t$ corresponding to $\bs u$ is assumed to be of the form 
\begin{align}
	P_{\mathrm{int}}(\cl P_t, \bs u) := \int_{\cl P_t} (\bs T \cdot \nabla \bs u + \bbG \cdot \nabla\nabla \bs u ) \, dv, \label{eq:1}
\end{align}
where $\bs T$ is a second order tensor (the \emph{stress}) and $\bbG$ is a third order tensor (the \emph{hyperstress}). It is assumed that $\bbG^{T(2,3)} = \bbG$, reflecting the symmetry of the second gradient of the test velocity. The internal power being frame indifferent dictates that $\bs T$ is symmetric. The \textit{virtual power expended by the external forces} on $\cl P_t$ corresponding to $\bs u$ is expressed as  
\begin{align}
	P_{\mathrm{ext}}(\cl P_t, \bs u) := \int_{\p \cl P_t} (\bs t \cdot \bs u + \bs m \cdot \p_n \bs u ) \, da + 
	\int_{\cl P_t} \rho (\bs b - \dot{\bs v})  \cdot \bs u \, dv, \label{eq:2}
\end{align}
where $\bs t$ is the \textit{traction} and $\bs m$ is the \textit{hypertraction} acting on $\p \cl P_t$. The second term in \eqref{eq:2} includes the specific inertial body force $\bs b^{in}$ given by the classical expression 
\begin{align}
	\bs b^{in} := -\dot{\bs v}. \label{eq:classicalin}
\end{align}
The principle of virtual work \cite{Germain73a, Germain73b, FriedGurtin06} states that for all $\cl P_t$ and test velocities $\bs u$ on $\cl B_t$, we have  
\begin{align}
	P_{\mathrm{int}}(\cl P_t, \bs u) = P_{\mathrm{ext}}(\cl P_t, \bs u). \label{eq:3}
\end{align}

From \eqref{eq:3}, integration by parts, and the fundamental lemma of the calculus of variations, \cite{FriedGurtin06} derives local equations expressing the local form of balance of forces (including inertial) and expressions for $\bs t$ and $\bs m$ in terms of $\bs T$, $\bbG$, and the geometry of $\p \cl P_t$. The local form of balance of forces is given by 
\begin{align}
	\rho \dot{\bs v} = \mathrm{div}\,(\bs T - \mathrm{div}\,\bbG) + \rho \bs b. \label{eq:6}
\end{align}  
The traction and hypertraction are then given by 
\begin{align}
	\bs t &= \bs T \bs n - (\mathrm{\div}\, \bbG) \bs n - \mathrm{div}_s(\bbG \bs n) - 2 K\bbG[\bs n \tens \bs n], \label{eq:4} \\
	\bs m &= \bbG[\bs n \tens \bs n]. \label{eq:5}
\end{align}
Above, $K$ is the \textit{mean curvature} of $\p \cl P_t$ and $\mathrm{div}_s$ is the \textit{surface divergence} on $\p \cl P_t$ (see, e.g., Chapter 2 of  \cite{steigmann2023lecture} for their expressions in general curvilinear coordinates). In the case that $\cl P_t = \cl B_t$, specifying $\bs t$ and/or $\bs m$ on parts of the boundary constitutes prescribing the tractions and hypertractions for the boundary value problem of interest. See Section 4 of \cite{FriedGurtin06} for more details. 

The work \cite{FriedGurtin06} also incorporates an additional spatial gradient of the velocity in the kinetic energy by appealing to a nonstandard \textit{principle of inertial power balance}. The form of virtual power expended by the external forces on a part $\cl P_t$ corresponding to a test velocity $\bs u$ is now expressed as
\begin{align}
	P_{\mathrm{ext}}(\cl P_t, \bs u) := \int_{\p \cl P_t} \Bigl (\bs t \cdot \bs u + \bs m \cdot \frac{\p \bs u}{\p \bs n} \Bigr ) \, da + 
\int_{\cl P_t} \rho \Bigl (\bs b - \dot{\bs v} + \ell_0^2 \mathrm{\div}\bigl (\dot{\overline{\grad\bs v}} \bigr ) \Bigr ) \cdot \bs u \, dv, \label{eq:7}
\end{align}
where $\ell_0 \geq 0$ is an intrinsic length scale. 
When compared to \eqref{eq:2}, we see that the specific inertial body force is now given by
\begin{align}
	\bs b^{in} = - \Bigl ( \dot{\bs v} - \ell_0^2 \mathrm{\div}\bigl (\dot{\overline{\grad\bs v}} \bigr ) \Bigr ). \label{eq:8}
\end{align} 
The local form of balance of forces then takes the form 
\begin{align}
	\rho \Bigl ( \dot{\bs v} - \ell_0^2 \mathrm{\div}\bigl (\dot{\overline{\grad\bs v}} \bigr ) \Bigr ) = \mathrm{div}\,(\bs T - \mathrm{div}\,\bbG) + \rho \bs b. \label{eq:9}
\end{align}
For $\cl P_t = \cl B_t$, one instead specifies on parts of the boundary the non-inertial traction $\bs t^{ni}$ and/or hypertraction $\bs m^{ni}$ for the boundary value problem of interest:
\begin{align}
\bs t^{ni} &:= \bs t + \rho \ell_0^2 \bigl (\dot{\overline{\grad\bs v}} \bigr ) \bs n, \label{eq:10} \\
\bs m^{ni} &:= \bs m, \label{eq:10b}
\end{align}
where $\bs t$ and $\bs m$ are given in \eqref{eq:4} and \eqref{eq:5}. See Section 12 of \cite{FriedGurtin06} for more details. 

\subsection{Non-simple, incompressible, linear viscous fluid models} Motivated by modeling incompressible viscous flows at small-length scales, the works \cite{FriedGurtin06, musesti2009isotropic} propose a generalization of the constitutive relations for the classical incompressible Navier-Stokes model of viscous fluids. In Cartesian components, they take the form   
\begin{align}
T_{ij} &= -p \delta_{ij} + T_{0ij}, \quad T_{0kk} = 0, \label{eq:11} \\
G_{ijk} &= -\frac{1}{2}\bigl (
\delta_{ij}\pi_k + \delta_{ik}\pi_j
\bigr ) + G_{0ijk}, \quad G_{0ijk} = G_{0ikj}, \quad G_{0iik} = 0. \label{eq:12} 
\end{align} 
Above, $p$ is the \textit{mean normal stress} and the vector $\bs \pi$ with Cartesian components $\{\pi_k\}_{k=1}^3$ is the so-called \textit{hyperpressure}. Following standard nomenclature and \cite{FriedGurtin06, musesti2009isotropic}, we will refer to $p$ as the \textit{pressure}; see \cite{rajagopal2015remarks} for an in-depth discussion of various notions of pressure. The tensor $\bs T_0$ with Cartesian components $\{T_{0ij}\}_{i,j=1}^3$ is the \textit{extra stress} and the tensor $\bbG_0$ with components $\{G_{0ijk}\}_{i,j,k=1}^3$ is the \textit{extra hyperstress}.  

The works \cite{FriedGurtin06, musesti2009isotropic} express $\bs T_0$ and $\bbG_0$ as linear isotropic functions of $\bs D$ and $\grad \grad \bs v$ respectively that are consistent with the constraints in \eqref{eq:11} and \eqref{eq:12}. These take the most general form \cite{musesti2009isotropic}:  
\begin{align}
	T_{0ij} &= 2\mu D_{ij}, \label{eq:extrastresscomponents} \\
  G_{0ijk} &= \eta_1 v_{i,jk} + \eta_2(v_{k,ij} + v_{j,ik} - v_{i,rr} \delta_{jk}) \\ &\quad+ \eta_3(v_{k,rr}\delta_{ij} + v_{j,rr}\delta_{ik} - 4v_{i,rr}\delta_{jk}). \label{eq:extrahyperstresscomponents}
\end{align}
In direct notation, the above reads
\begin{align}
\bs T_0 &= 2\mu \bs D, \\
\bbG_0 &= \eta_1 \nabla \nabla \bs v + \eta_2 \Bigl [
\nabla \nabla \bs v^{T(1,3)} + \nabla \nabla \bs v^{T(1,2)} - \Delta \bs v \tens \bs I  
\Bigr ] \\
&\quad + \eta_3 \Bigl [
\bs I \tens \Delta \bs v + (\bs I \tens \Delta \bs v)^{T(2,3)} - 4 \Delta \bs v \tens \bs I 
\Bigr ]. \label{eq:hyperstressdirect}
\end{align}
Using the fact that 
\begin{align}
	G_{0ijk} &= (\eta_1 + 2 \eta_2)D_{ij,k} + (\eta_1 - 2\eta_2)W_{ij,k} - (\eta_2 + 4 \eta_3)v_{i,rr}\delta_{jk} \\&\quad+ 
	\eta_3(v_{k,rr} \delta_{ij} + {\delta_{ik}} v_{j,rr}), 
\end{align}
The dissipation rate $\xi$ is given by
\begin{align}
	\xi &:= \bs T \cdot \bs D + \bbG\cdot \nabla \nabla \bs v \\
	&= 2\mu |\bs D|^2 + (\eta_1 + 2 \eta_2)|\nabla \bs D|^2 + (\eta_1 - 2 \eta_2)|\nabla \bs W|^2 - (\eta_2 + 4 \eta_3)|\Delta \bs v|^2. \label{eq:dissrate}
\end{align}
We split $\nabla \bs D$ and $\nabla \bs W$ into their deviatoric and spherical parts with respect to the first and third indices. We do this while also preserving their symmetries in the first and second indices to maintain independence of their components.  More precisely, we write: 
\begin{gather}
	\nabla \bs D = \widehat{\nabla \bs D} + \nabla \bs D_{s}, \quad \nabla \bs W = \widehat{\nabla \bs W} + \nabla \bs W_{s}, \\
	(\widehat{\nabla \bs D})_{ijk} := D_{ij,k} - \frac{1}{4}(\delta_{ik}D_{\ell j, \ell} + \delta_{jk}D_{\ell i,\ell}), \\ 
	(\widehat{\nabla \bs W})_{ijk} := W_{ij,k} - \frac{1}{2}(\delta_{ik}W_{\ell j, \ell} - \delta_{jk}W_{\ell i,\ell}). 
\end{gather}
Using that $\mathrm{div}\, \bs v = 0$, it follows that 
\begin{align}
	D_{\ell j, \ell} = - W_{\ell j, \ell} = \frac{1}{2} \Delta v_j. 
\end{align}
Thus, the dissipation rate is given as a sum of three independent terms,
\begin{align}
	\xi &= 2\mu|\bs D|^2 + (\eta_1 + 2 \eta_2) |\widehat{\nabla \bs D}|^2 + (\eta_1 - 2 \eta_2) |\widehat{\nabla \bs W}|^2 \\&\quad+
	\frac{1}{8}(3 \eta_1 - 10\eta_2 - 32 \eta_3)|\Delta \bs v|^2. 
\end{align}
Then $\xi \geq 0$ for all values of $\nabla \nabla \bs v$ with $\mathrm{div}\, \bs v = 0$ if and only if 
\begin{align}
\eta_1 \geq 2 |\eta_2|, \quad \frac{1}{8}(3 \eta_1 - 10\eta_2 - 32 \eta_3) \geq 0. \label{eq:constantconstraints}	
\end{align}

We define the following independent length scales of the model,
\begin{gather}
	\ell_2 := \Bigl [\frac{1}{2\mu}(\eta_1 + 2 \eta_2) \Bigr ]^{1/2}, \quad \ell_3 := \Bigl [ \frac{1}{2\mu}(\eta_1 - 2 \eta_2)\Bigr ]^{1/2}, \\ \ell_4 := \Bigl [ \frac{1}{16\mu}(3 \eta_1 - 10\eta_2 - 32 \eta_3)\Bigr ]^{1/2}, \label{eq:lengthscales}
\end{gather}
and we note that 
\begin{align}
\eta_1 = \mu(\ell_2^2 + \ell_3^2), \quad \eta_2 = \frac{1}{2}\mu(\ell_2^2-\ell_3^2), \\
\eta_3 = \mu \Bigl (
-\frac{1}{16} \ell_2^2 + \frac{1}{4}\ell_3^2 - \frac{1}{2} \ell_4^2
\Bigr ). 
\end{align}
The dissipation rate then takes the form
\begin{align}
	\xi = 2\mu |\bs D|^2 + 2\mu \ell_2^2 |\widehat{\nabla \bs D}|^2 + 2\mu \ell_3^2 |\widehat{\nabla \bs W}|^2 +
2 \mu \ell_4^2 |\Delta \bs v|^2. 	
\end{align}
Motivated by the ultimate form of the governing field equations (see \eqref{eq:field}), we define a primary length scale $\ell_1$ via 
\begin{align}
 \ell_1^2 &:= \frac{1}{\mu}(\eta_1 - \eta_2 - 4 \eta_3)\\
	&= \frac{1}{2\mu}\eta_1 + \frac{1}{8\mu}(\eta_1 + 2\eta_2) + \frac{1}{8\mu}(3 \eta_1 - 10\eta_2 - 32 \eta_3) \\
	&= \frac{3}{4}\ell_2^2 + \frac{1}{2}\ell_3^2 + 2\ell_4^2 \label{eq:primelength}
\end{align}
We note that
\begin{align}
	\ell_1 = 0 \iff \ell_2 = \ell_3 = \ell_4 = 0, \label{eq:ell10}
\end{align}
In the mathematically simplifying case that the dissipation rate \textit{depends only on the spherical part} of the second gradient of the velocity, there is only one length scale appearing in the dissipation rate, 
\begin{align}
	\ell_2 = \ell_3 = 0, \quad \ell_1^2 = 2\ell_4^2, 
\end{align}
and
\begin{align}
	\xi = 2\mu |\bs D|^2 + 2\mu \ell_1^2 |\Delta \bs v|^2. 
\end{align}


\subsection{The proposed model: specifying the hyperstress}
The theory presented in the previous section raises several issues that must be addressed {in order} to develop a well-defined and predictive model of incompressible viscous flows incorporating small-length scale effects. The first challenge, common in higher gradient continuum models (see, e.g., \cite{Askes2011}), is the determination of the additional length scales $\ell_2, \ell_3, \ell_4$. The second, specific to the present model, concerns the hyperstress vector $\bs \pi$, which remains indeterminate from the governing equations. Indeed, substituting \eqref{eq:11}, \eqref{eq:12}, \eqref{eq:extrastresscomponents}, and \eqref{eq:extrahyperstresscomponents} into \eqref{eq:9} yields  
\begin{align}
	\rho \Bigl ( \dot{\bs v} - \ell_0^2 \mathrm{\div}\bigl (\dot{\overline{\grad\bs v}} \bigr ) \Bigr ) =  - \nabla(p - \mathrm{div}\, \bs \pi) + \mu \Delta \bs v - \mu \ell_1^2 \Delta^2 \bs v + \rho \bs b, \label{eq:13}
\end{align}  
where $\ell_1$ is given in \eqref{eq:primelength}. Here, only the quantity $p - \mathrm{div}\,\bs \pi$, known as the \emph{effective pressure} \cite{FriedGurtin06}, appears in the equations of motion and can be determined (given appropriate boundary conditions) by taking the divergence of \eqref{eq:13}. However, knowledge of $\bs \pi$ is necessary for boundary conditions posed in terms of the hypertraction \eqref{eq:5} or \eqref{eq:10b}. Thus, a constitutive relation for $\bs \pi$ must be prescribed.

The work \cite{krawietz2022surface} gives a compelling argument connecting $\bs \pi$ and $\grad p$. More precisely, by taking a certain incompressible limit\footnote{More precisely, the absence of a coupling term in a compressible second-gradient fluid's quadratic free energy and the vanishing of two retardation times in the incompressible limit lead to \eqref{eq:hyperpress}. See \cite{krawietz2022surface} for more details.}, it follows that 
\begin{align}
	\bs \pi = \ell^2 \grad p, \label{eq:hyperpress}
\end{align}   
where $\ell^2$ is an additional length scale of the problem; see Section 4 of \cite{krawietz2022surface}. Since the model presented already contains four independent length scales ($\ell_0, \ell_2, \ell_3, \ell_4$), we will make the mathematical simplifying choice that $\ell = \ell_1$ where $\ell_1$ is given in \eqref{eq:primelength}. Then the governing field equations become 
\begin{equation}
\rho \Bigl ( \dot{\bs v} - \ell_0^2 \mathrm{\div}\bigl (\dot{\overline{\grad\bs v}} \bigr ) \Bigr ) = - \nabla(p - \ell_1^2 \Delta p) + \mu \Delta \bs v - \mu\ell_1^2\Delta^2 \bs v + \rho \bs b. \label{eq:field}
\end{equation} 

In summary, the full model that we propose in this work contains four independent length scales with
\begin{align}
\bs T &= -p \bs I + 2\mu \bs D, \\
\bbG &= -\frac{\ell_1^2}{2}\Bigl [\bs I \tens \nabla p + (\bs I \tens \nabla p)^{T(2,3)} \Bigr ] + \mu(\ell_2^2 + \ell_3^2) \nabla \nabla \bs v \\&\quad+ \frac{1}{2}\mu(\ell_2^2 -\ell_3^2) \Bigl [
\nabla \nabla \bs v^{T(1,2)} + \nabla \nabla \bs v^{T(1,3)} - \Delta \bs v \tens \bs I  
\Bigr ] \\
&\quad +  \mu\Bigl (-\frac{1}{16}\ell_2^2 + \frac{1}{4} \ell_3^2 - \frac{1}{2} \ell_4^2 \Bigr )\Bigl [
\bs I \tens \Delta \bs v + (\bs I \tens \Delta \bs v)^{T(2,3)} - 4 \Delta \bs v \tens \bs I 
\Bigr ], \label{eq:Gequations}
\end{align}
where 
\begin{align}
	\ell_1^2 = \frac{3}{4}\ell_2^2 + \frac{1}{2}\ell_3^2 + 2\ell_4^2.
\end{align}
In contrast to the isotropic model derived in \cite{krawietz2022surface}, \eqref{eq:Gequations} \textit{contains two fewer free constants} due to: splitting the hyperstress into {an expression containing} the hyperpressure and {the} extra hyperstress at the outset, and setting the length scale appearing in \eqref{eq:hyperpress} to be the same as $\ell_1$. We note that we recover the classical Navier-Stokes model when $\ell_0 = \ell_1 = 0$ by \eqref{eq:ell10} and \eqref{eq:field}. 

Finally, we remark that in the simplifying case when the extra hyperstress \textit{depends only on the spherical part} of the second gradient of the velocity, $\bbG$ takes the mathematically attractive form
\begin{align}
	\bbG &= -\frac{\ell_1^2}{2}\Bigl [\bs I \tens \nabla p + (\bs I \tens \nabla p)^{T(2,3)} \Bigr ] + 2 \mu \ell_4^2 \Delta \bs v \tens \bs I \\
	&= -\frac{\ell_1^2}{2}\Bigl [\bs I \tens \nabla p + (\bs I \tens \nabla p)^{T(2,3)} \Bigr ] + \mu \ell_1^2 \Delta \bs v \tens \bs I, \label{eq:spherical}
\end{align}
and contains only one unknown length scale.

\subsection{Equation governing the pressure}
Taking the divergence of both sides of \eqref{eq:field}  and using the incompressibility assumption implies that the pressure satisfies
\begin{equation}
	\rho \mathrm{div} \Bigl ( \dot{\bs v} - \ell_0^2 \mathrm{\div}\bigl (\dot{\overline{\grad\bs v}} \bigr ) \Bigr ) = -\Delta p + \mu\ell_1^2\Delta^2 p + \rho \mathrm{div}\, \bs b.
\end{equation}
Simplifying the left hand side of the previous equation gives
\begin{equation}
	\rho \Bigl (\nabla \bs v \cdot \nabla \bs v^T - \ell_0^2 \nabla \nabla \bs v \cdot \nabla \nabla \bs v^{T(1,2)} 
	- \ell_0^2 \nabla \Delta \bs v \cdot \nabla \bs v^T \Bigr ) = -\Delta p + \mu\ell_1^2\Delta^2 p + \rho \mathrm{div}\, \bs b. \label{eq:pressureequation}
\end{equation}

\subsection{Boundary conditions}\label{s:boundary}
The work \cite{FriedGurtin06} discusses a broad range of both no-slip and slip conditions on a fixed smooth piece of the boundary $\cl S \subseteq \p \cl B_t$ for a so-called \textit{passive environment}. In what follows, we are assuming that $\ell_1 > 0$. For no-slip conditions, \cite{FriedGurtin06} gives a one parameter family of generalized adherence conditions parameterized by a length scale $\ell > 0$: 
\begin{align}
	\bs v = \bs 0, \quad \bs m = -\mu \ell \frac{\p \bs v}{\p \bs n}, \quad \mbox{on }\cl S,   
\end{align}
where $\bs m$ is given in \eqref{eq:5}. The two extreme limits $\ell \rar 0$ and $\ell \rar \infty$ correspond to what \cite{FriedGurtin06} refers to as \textit{weak adherence} and \textit{strong adherence}, given respectively on $\cl S$ by 
\begin{gather}
	\bs v = \bs 0, \quad \bs m = \bs 0 \quad \mbox{(weak adherence)}, \\
	\bs v = \bs 0, \quad \frac{\p \bs v}{\p \bs n} = \bs 0 \quad \mbox{(strong adherence)}.
\end{gather}
The study \cite{FriedGurtin06} also generalizes the classical Navier slip condition $\bs P \bs T \bs n = -\frac{\mu}{\lambda}\bs v$ to 
\begin{align}
	\bs P \bs t = -\frac{\mu}{\lambda}\bs v, \quad \bs m = \bs 0, \quad \quad \mbox{on }\cl S,
\end{align}  
where $\bs P$ is the projection onto to the tangent space of $\cl S$, $\bs t$ is given in \eqref{eq:4}, and $\lambda$ is the \textit{slip length}. In this work, we will focus on weak adherence and strong adherence boundary conditions.

We observe that for our model, we must pose additional boundary conditions to determine the pressure uniquely on a domain with a fixed boundary $\cl S$. Indeed, \eqref{eq:pressureequation} is a fourth-order elliptic equation rather than the classical second order elliptic equation corresponding to $\ell_0 = \ell_1 = 0$. We propose that in addition to prescribing $p$ on $\cl S$, we assume that 
\begin{align}
	\frac{\p p}{\p \bs n} = 0, \quad \mbox{on }\cl S. \label{eq:pressureboundary}
\end{align}  
This assumption is motivated by \eqref{eq:5}, \eqref{eq:Gequations}, and the fact that $\bs m = \bs 0$ appears in the case of both weak adherence and the generalized slip conditions.

\section{A second-gradient model for incompressible viscous fluids with pressure dependent viscosities} 

\subsection{An earlier model} Previous studies have examined a modification of the incompressible Navier–Stokes model, rooted in \textit{classical continuum mechanics}, that incorporates pressure-dependent viscosity in an effort to address the limitations of the classical theory in high-pressure regimes:
\begin{align}
	\rho \dot{\boldsymbol{v}} = \div \boldsymbol{T} + \rho \bs b, \quad \boldsymbol{T} = -p\boldsymbol{I} + 2\hat{\mu}(p) \boldsymbol{D}, \quad \div \bs v = 0. \label{eq:pressuredependent}
\end{align}
A commonly used relation between viscosity and pressure that is consistent with experiment is {Barus' exponential formula} \cite{Barus1893}: $\hat \mu(p) = \mu_0 \exp(\alpha(p - p_0))$ where $\mu_0 > 0$, $\al > 0$, and $p_0$ are fixed constants.  Several works have investigated the modeling and simulation aspects of \eqref{eq:pressuredependent} \cite{Bair1998, Rajagopal2003, Bayada2013, Gustafsson2015, Denn1981, Hron2001, Vasudevaiah2005, Prusa2010, Kalogirou2011, Rajagopal2012, Rajagopal2013, Marusic-Paloka2013, Housiadas2015, PrusaRehor2016, Gwynllyw1996, Hron2003, Chung2010, Knauf2013, Janecka2014}. However, the well-posedness theory of \eqref{eq:pressuredependent} is highly restricted due to \textit{the potential for the equation to change type} (see \eqref{eq:pressureequationnew}). Local well-posedness is ensured under the condition that $\hat{\mu}(p)/p \to 0$ as $p \to \infty$ which contradicts Barus' exponential formula and other experimental findings \cite{Renardy1986}. Additionally, it is only known that sufficiently slow steady states exist and are unique within a fixed bounded domain \cite{GazzolaSecchi1998}. The narrow scope of these mathematical results and the potential for the equation to switch type raise serious questions about the general predictability of the model \eqref{eq:pressuredependent}.

\subsection{A second-gradient model with pressure dependent viscosity}
We propose a second-gradient model for incompressible viscous fluids with pressure dependent viscosities that is governed by \eqref{eq:9} with 
\begin{align}
	\bs T &= -p \bs I + 2\hat \mu(p) \bs D, \label{eq:pdT}\\
	\bbG &= -\frac{\ell_1^2}{2}\Bigl [\bs I \tens \nabla p + (\bs I \tens \nabla p)^{T(2,3)} \Bigr ] + \mu_0(\ell_2^2 + \ell_3^2) \nabla \nabla \bs v \\&\quad+ \frac{1}{2}\mu_0(\ell_2^2 -\ell_3^2) \Bigl [
	\nabla \nabla \bs v^{T(1,2)} + \nabla \nabla \bs v^{T(1,3)} - \Delta \bs v \tens \bs I  
	\Bigr ] \\
	&\quad +  \mu_0\Bigl (-\frac{1}{16}\ell_2^2 + \frac{1}{4} \ell_3^2 - \frac{1}{2} \ell_4^2 \Bigr )\Bigl [
	\bs I \tens \Delta \bs v + (\bs I \tens \Delta \bs v)^{T(2,3)} - 4 \Delta \bs v \tens \bs I 
	\Bigr ], \label{eq:pressureGequations}
\end{align} 
with $\mu_0 = \hat \mu(p_0)$ for a fixed value of pressure $p_0$ and $\ell_1 > 0$. The governing field equations then read
\begin{align}
	\rho \bigl [\dot{\boldsymbol{\bs v}} - \ell_0^2 \div(\dot{\overline{\nabla \bs v}}) \bigr ] = -\grad( p-\ell_1^2 \Delta p) + \div ( 2 \hat \mu(p) \bs D) - \ell_1^2 \mu_0 \Delta^2 \bs v + \rho \bs b, \quad \label{eq:Gpressuredependent} 
\end{align}
with $\mathrm{div}\, \bs v = 0$. We observe that for a simple shear flow $\bs v = \gamma x_2 \bs e_1$ under constant pressure $p_1$, we have $\bbG = \bs 0$. If $\p \cl P_t \subseteq \{ x_2 = \mathrm{constant} \}$, then on $\cl P_t$ with $\bs n = -\bs e_2$, $\bs m = \bs 0$ and 
\begin{align}
	\bs t = \bs T (-\bs e_2) = -\hat \mu(p_1) \gamma \bs e_1 + p_1 \bs e_2. \label{eq:stressshear}  
\end{align}
Thus, the shear stress exerted by the (bottom) layer of fluid adjacent to the (top) layer of the fluid along $\p \cl P_t$ is given by $-\hat \mu(p_1) \gamma \bs e_1$. This is the same form as what would be obtained from the classical Navier-Stokes model but with a generalized viscosity depending on the given constant density $p_1$. It is in this sense that we say \eqref{eq:pdT} and \eqref{eq:pressureGequations} model a fluid with pressure dependent viscosity. 

The equation determining the pressure is obtained by taking the divergence of $\eqref{eq:Gpressuredependent}$, yielding 
\begin{gather}
\rho \Bigl (\nabla \bs v \cdot \nabla \bs v^T - \ell_0^2 \nabla \nabla \bs v \cdot \nabla \nabla \bs v^{T(1,2)} 
- \ell_0^2 \nabla \Delta \bs v \cdot \nabla \bs v^T \Bigr )\\ = \ell_1^2 \Delta^2 p - \div \bigl [ (\bs I - 2\hat \mu'(p) \bs D) \grad p \bigr ] + \hat \mu'(p) \nabla p \cdot \Delta \bs v + \rho \div \bs b. \label{eq:pressureequationnew}
\end{gather}
Due to the term $\ell_1^2 \Delta^2 p$ appearing in \eqref{eq:pressureequationnew}, it follows that equation \eqref{eq:pressureequationnew} is \emph{always elliptic in $p$} in contrast with the model \eqref{eq:pressuredependent} ($\ell_0 = \ell_1 = 0$) which requires definiteness of $\bs I - 2\hat \mu'(p) \bs D$ for ellipticity. This distinction underscores the potential predictive advantage of \eqref{eq:Gpressuredependent} over \eqref{eq:pressuredependent}. The remainder of this work will focus on the constant viscosity case, and the model defined by \eqref{eq:pressuredependent}, \eqref{eq:pressureGequations}, and \eqref{eq:Gpressuredependent} will be explored in forthcoming work.

\section{Poiseuille Flow}

In this section, we consider the classical problem for Poiseuille flow for our model: steady unidirectional flow in a cylindrical tube of radius $R$ induced by a constant pressure gradient. We work in cylindrical coordinates $(r,\theta,z)$ with natural basis vectors $\{\bs g_i \}$, dual basis vectors $\{\bs g^i \}$, and orthonormal basis vectors $\bs e_i = \bs g_i/|\bs g_i| = \bs g^i/|\bs g^i|$ with $i \in \{r, \theta, z\}$. We assume that the velocity $\bs v$ and pressure gradient $\grad p$ take the forms
\begin{align}
	\bs v = v(r)\bs e_z, \quad \nabla \bs p = -\beta \bs e_z. \label{eq:14}
\end{align}
Then 
\begin{gather}
	\nabla \bs v = \bs v_{,i} \tens \bs g^i
	= v' \bs e_z \tens \bs e_r, \\
	\nabla \nabla \bs v = (\nabla \bs v)_{,i}\tens \bs g^i
	= v'' \bs e_z \tens \bs e_r \tens \bs e_r + \frac{1}{r}v' \bs e_z \tens \bs e_\theta \tens \bs e_\theta, 
\end{gather}
where $\mbox{}' := \frac{d}{dr}$. In this section and the next, we assume that $\ell_1 > 0$. 

\subsection{General solution}
The field equations \eqref{eq:field} reduce to
\begin{equation}
    -\frac{\beta}{\mu} = \calL(1 - \ell_1^2 \calL)v.
\end{equation}
where $\cl L = \frac{d^2}{dr^2} + \frac{1}{r} \frac{d}{dr}$ and $\ell_1$ is defined given by \eqref{eq:primelength}. 
The above equation has the following general solution:
\begin{equation}
  v = \frac{\beta}{4\mu}(R^2 - r^2) + c_1\ell_1^2 I_0(r / \ell_1) + c_2 \ell_1^2 K_0(r / \ell_1) + c_3\log r + c_4,
\end{equation}
where the first term is the classical solution for Poiseuille flow with no-slip conditions, and $I_0$ and $K_0$ are modified Bessel functions of the first kind {of order $0$}. In order for the solution to be regular at $r = 0$, we must have that $c_2 = c_3 = 0$.

Thus, after relabeling the constants, the velocity field is give{n} by
\begin{equation}
    v(r)  = \frac{\beta}{4\mu}(R^2 - r^2) +  c_1\ell_1^2 I_0(r / \ell_1) + c_2
\end{equation}
Since both weak adherence and strong adherence boundary conditions include the no-slip condition, we must have that 
\begin{equation} \label{eqn:PNoSlip}
    0 = v(R) =  c_1\ell_1^2 I_0(R/\ell_1) + c_2,
\end{equation} 
and thus,
\begin{align}
	v(r) = \frac{\beta}{4\mu}(R^2 - r^2) + c_1 \ell_1^2 \Bigl [ I_0(r/\ell_1) - I_0(R/\ell_1) \Bigr ]. \label{eq:vequation}
\end{align}
\subsection{Strong adherence boundary conditions}
In the strong adherence case, we have that $\p_{\bs n} \bs v = \bs 0$ on the boundary of the cylinder. In our setting this is equivalent to $v'(R) = 0$. Thus,
\begin{equation*}
    v'(R) = c_1\ell_1 I_0'(R/\ell_1) - \frac{\beta}{2\mu}R = 0 \implies c_1 = \frac{\beta}{2\mu}\frac{R}{\ell_1 I_0'(R/\ell_1)}
\end{equation*}
Substituting $c_1$ into \eqref{eq:vequation} gives us the final expression for the velocity satisfying strong adherence boundary conditions:
\begin{equation}
    v(r) = \frac{\beta}{4\mu}(R^2 - r^2) + \frac{\beta}{2\mu}\frac{R\ell_1}{I_0'(R/\ell_1)}(I_0(r/\ell_1) - I_0(R/\ell_1)). \label{eq:strongadherencePois}
\end{equation}

The discharge rate is defined via 
\begin{align}
	Q := \iint_{A} \bs v \cdot \bs n \, dA
\end{align}
where $A$ is the cross-section of the pipe and $\bs{n}$ is the unit normal to $A$. We conclude from \eqref{eq:strongadherencePois} and the power series representation for the modified Bessel functions of the first kind \cite{abramowitz1964handbook},
\begin{align}
I_{\nu}(\rho) = \Bigl (\frac{1}{2}\rho \Bigr )^\nu \sum_{k = 0}^\infty\frac{(\rho^2/4)^k}{k!\Gamma(\nu+k+1)},
\end{align}
that
\begin{align}
	Q &= 2\pi \int_0^R \Bigl [ \frac{\beta}{4\mu}(R^2 - r^2) + \frac{\beta}{2\mu}\frac{R\ell_1}{I_0'(R/\ell_1)}(I_0(r/\ell_1) - I_0(R/\ell_1)) \Bigr ] r dr \\
	&= \frac{\pi \beta R^4}{8\mu} + \frac{\pi \beta}{\mu}\frac{R\ell_1}{I_0'(R/\ell_1)}\Bigl (R \ell_1 I_1(R/\ell_1) - \frac{1}{2} R^2 I_0(R/\ell_1) \Bigr ). 
\end{align}

The dimensionless velocity $u(\sigma) := \frac{4\mu}{\beta R^2} v(R\sigma)$ as a function of the dimensionless variable $\sigma := r/R$ then takes the form
\begin{align}
	u(\sigma) = 1 - \sigma^2 + \frac{2 \lambda_1}{I_0'(1/\lambda_1)}(I_0(\sigma/\lambda_1) - I_0(1/\lambda_1)), \label{eq:dimensionlessstrongP}
\end{align}
where $\lambda_1 := \ell_1/R$ is the nondimensionalized length scale. In addition, the dimensionless discharge rate $\Phi := \frac{8\mu Q}{\pi \beta R^4}$ is given by  
\begin{align}
	\Phi = 1 + \frac{8\la_1}{I_0'(1/\la_1)}\Bigl [ \la_1 I_1(1/\la_1) - \frac{1}{2} I_0(1/\la_1) \Bigr ].  
\end{align}

\begin{figure}[h]
	\centering
	\includegraphics[trim = {0 6cm 0 6cm},clip,scale = 0.55]{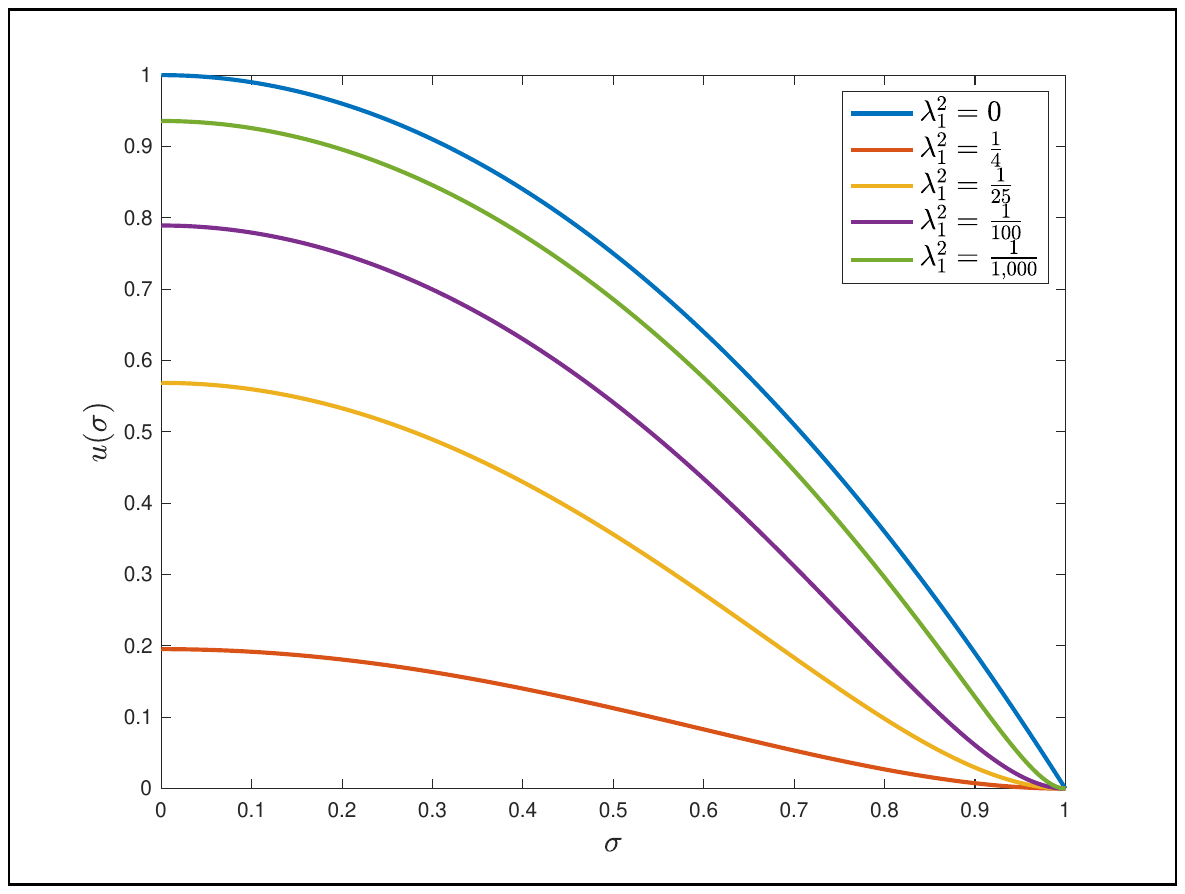}
	\caption{Graph of $u(\sigma)$ for Poisuelle flow with strong adherence boundary conditions. The curve for $\lambda_1 = 0$ corresponds to the classical solution.}
	\label{fig:f1}
\end{figure}

\subsection{Weak adherence boundary conditions}
In the weak adherence case, we assume that $\bbG[\bs n \tens \bs n] = \bs 0$. Since we require $p'(R) = 0$, the velocity field must satisfy
\begin{align}
	(\eta_1 - \eta_2 - 4 \eta_3)v''(R) - (\eta_2 + 4 \eta_3)\frac{1}{R}v'(R) = 0. \label{eq:weakadherencebdy1}
\end{align}
The boundary condition \eqref{eq:weakadherencebdy1} can be written in terms of the length scales $\{\ell_j \}_{j = 1}^4$ defined in Section 2.3 via
\begin{align}
	\mu \ell_1^2 v''(R) - \mu \Bigl (
	\frac{1}{4} \ell_2^2 + \frac{1}{2} \ell_3^2 - 2 \ell_4^2
	\Bigr )\frac{1}{R}v'(R) = 0,   
\end{align}
where $\ell_1^2 = \frac{3}{4}\ell_2^2 + \frac{1}{2}\ell_3^2 + 2 \ell_4^2.$
In this case, we see that
\begin{equation*}
	v''(R) = -\frac{\beta}{2\mu} + c_1 I_0''(R/\ell_1),
\end{equation*}
so we can solve for $c_1$ when substituting these values into the boundary condition:
\begin{align*}
	\mu \ell_1^2 \Bigl ( -\frac{\beta}{2\mu} &+ c_1 I_0''(R/\ell_1) \Bigr ) - \mu \Bigl (
	\frac{1}{4} \ell_2^2 + \frac{1}{2} \ell_3^2 - 2 \ell_4^2
	\Bigr )\frac{1}{R} \Bigl ( c_1\ell_1 I_0'(R/\ell_1) - \frac{\beta}{2\mu}R \Bigr ) = 0 \\
	\implies & c_1 = \frac{R \beta (\ell_1^2 - \frac{1}{4}\ell_2^2 - \frac{1}{2}\ell_3^2+ 2\ell_4^2)}{2\mu(\ell_1^2RI_0''(R/\ell_1) - \ell_1(\frac{1}{4}\ell_2^2 + \frac{1}{2}\ell_3^2 - 2\ell_4^2)I_0'(R/\ell_1))}
\end{align*}
Then substituting $c_1$ into \eqref{eq:vequation} gives the following expression for the velocity satisfying the weak adherence boundary conditions:
\begin{equation}
	v(r) = \frac{\beta}{4\mu}(R^2 - r^2) + \frac{R \beta (\ell_1^2 - \frac{1}{4}\ell_2^2 - \frac{1}{2}\ell_3^2 + 2\ell_4^2)\ell_1^2 \bigl [ I_0(r/\ell_1) - I_0(R/\ell_1) \bigr ]}{2\mu(\ell_1^2RI_0''(R/\ell_1) - \ell_1(\frac{1}{4}\ell_2^2 + \frac{1}{2}\ell_3^2 - 2\ell_4^2)I_0'(R/\ell_1))} 
\end{equation}

Nondimensionalizing as in the previous section, with $\lambda_i := \ell_i/R$, our dimensionless velocity $u(\sigma)$ takes the form
\begin{align}
	u(\sigma) = 1 - \sigma^2 + \frac{2(\lambda_1^2-\frac{1}{4}\lambda_2^2-\frac{1}{2}\lambda_3^2 + 2\lambda_4^2) \lambda_1^2\bigl [ I_0(\sigma/\lambda_1) - I_0(1/\lambda_1)\bigr ]}{\lambda_1^2I_0''(1/\lambda_1) - \lambda_1(\frac{1}{4}\lambda_2^2 + \frac{1}{2}\lambda_3^2 - 2\lambda_4^2)I_0'(1/\lambda_1)},
\end{align}
and the  dimensionless discharge rate 
\begin{align}
\Phi &:= \frac{8 \mu}{\pi \beta R^4} Q \\
&= 1 + \frac{8(\lambda_1^2-\frac{1}{4}\lambda_2^2-\frac{1}{2}\lambda_3^2 + 2\lambda_4^2)\lambda_1^2\bigl [ \la_1 I_1(1/\lambda_1) - \frac{1}{2} I_0(1/\lambda_1)\bigr ]}{\lambda_1^2I_0''(1/\lambda_1) - \lambda_1(\frac{1}{4}\lambda_2^2 + \frac{1}{2}\lambda_3^2 - 2\lambda_4^2)I_0'(1/\lambda_1)}.
\end{align}

\begin{figure}[t]
	\centering
	\includegraphics[trim = {0 6cm 0 6cm},clip,scale = 0.55]{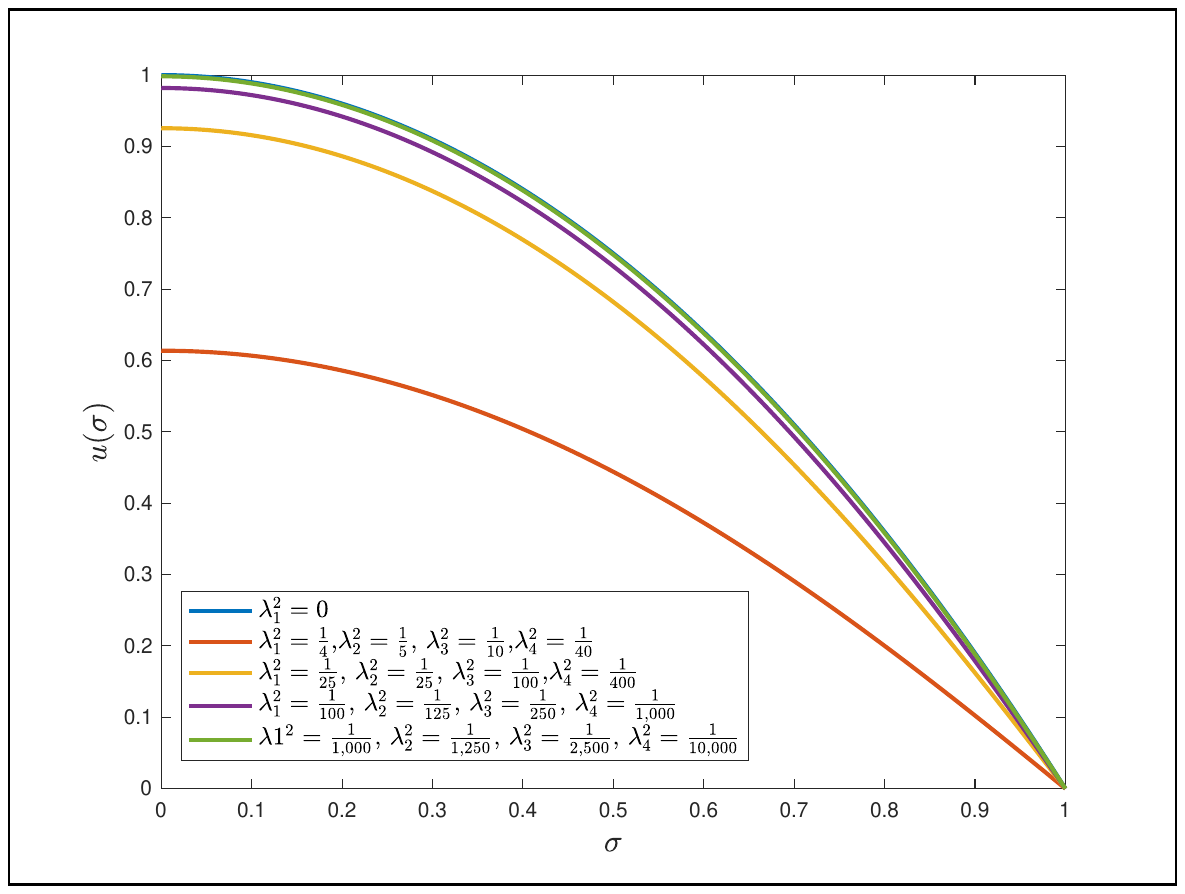}
	\caption{Graph of $u(\sigma)$ for Poisuelle flow with weak adherence boundary conditions. Here, $\lambda_2, \lambda_3 \neq 0$.}
	\label{fig:f2}
\end{figure}

\begin{figure}[t]
	\centering
	\includegraphics[trim = {0 6cm 0 6cm},clip,scale = 0.55]{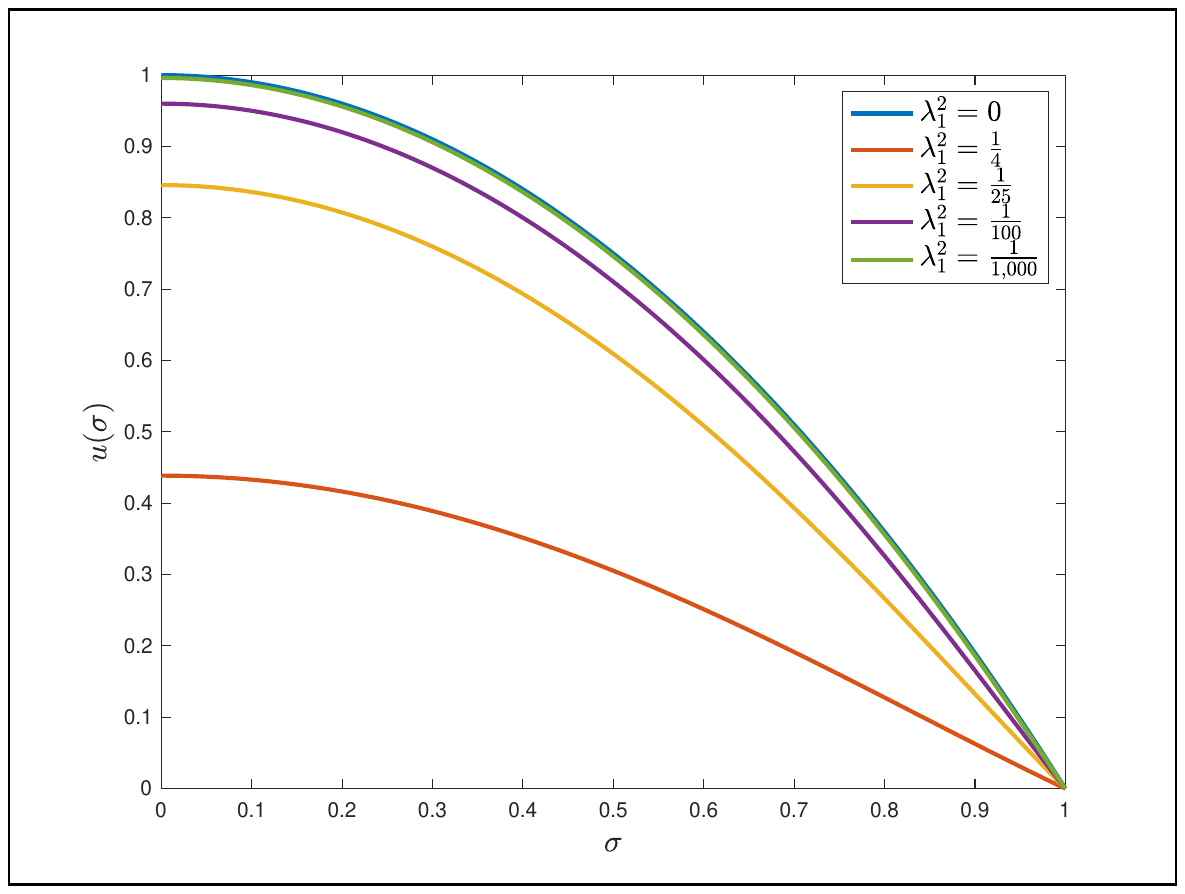}
	\caption{Graph of $u(\sigma)$ for Poisuelle flow {with weak adherence boundary conditions}. Here, $\lambda_2 = \lambda_3 = 0$.}
	\label{fig:f3}
\end{figure}

\subsection{Convergence of the velocity field to the solution for the classical Navier-Stokes model}

Recall our nondimensionalized solution in the strong adherence case:
\begin{align}
	u(\sigma) = 1 - \sigma^2 + \frac{2 \lambda_1}{I_0'(1/\lambda_1)}(I_0(\sigma/\lambda_1) - I_0(1/\lambda_1)). \label{eq:pvelocity}
\end{align}
Motivated by Figure \ref{fig:f1}, we will now show pointwise convergence of \eqref{eq:pvelocity} to the classical Navier-Stokes solution, $u_0(\sigma) = 1-\sigma^2$. Clearly, $u(1) = u_0(1) = 0$ so we will consider $\sigma \in (0,1)$. The argument for $\sigma = 0$ is similar. We will use the following large $z$ asymptotics for $I_0(z)$ \cite{abramowitz1964handbook}:
\begin{equation*}
	I_0(z) = \frac{e^z}{(2\pi z)^{1/2}}(1+O(z^{-1})) , \qquad I_0'(z) = \frac{e^z}{(2\pi z)^{1/2}}(1+O(z^{-1})).
\end{equation*}
We then can conclude that
\begin{align*}
	&\frac{2 \lambda_1 (I_0(\sigma/\lambda_1) - I_0(1/\lambda_1))}{I_0'(1/\lambda_1)} \\
    &= 2 \lambda_1 \left[ \frac{e^{\sigma/\lambda_1} \lambda_1^{1/2} \sigma^{-1/2}(1+O(\lambda_1/\sigma)) - e^{1/\lambda_1}\lambda_1^{1/2}(1 + O(\lambda_1))}{e^{1/\lambda_1}\lambda_1^{1/2} (1 + O(\lambda_1))}\right] \\
	&= 2 \lambda_1 \left[\frac{e^{-(1 - \sigma)/\lambda_1}\sigma^{-1/2}(1+O(\lambda_1/\sigma)) - (1 + O(\lambda_1))}{1+O(\lambda_1)}\right] \to 0 \text{ as } \lambda_1 \to 0.
\end{align*}
From the above calculation, we see that 
\begin{equation*}
    u(\sigma) = 1 - \sigma^2 + \frac{2 \lambda_1}{I_0'(1/\lambda_1)}(I_0(\sigma/\lambda_1) - I_0(1/\lambda_1)) \to 1 - \sigma^2 = u_0(\sigma) \text{ as } \lambda_1 \to 0,
\end{equation*}
proving pointwise convergence to the classical solution.

We now consider setting of weak adherence boundary conditions. Recall our nondimensionalized solution in the weak adherence setting:
\begin{align}
	u(\sigma) = 1 - \sigma^2 + \frac{2\lambda_1^2(\lambda_1^2-\frac{1}{4}\lambda_2^2-\frac{1}{2}\lambda_3^2 + 2\lambda_4^2)\bigl [ I_0(\sigma/\lambda_1) - I_0(1/\lambda_1)\bigr ]}{\lambda_1^2I_0''(1/\lambda_1) - \lambda_1(\frac{1}{4}\lambda_2^2 + \frac{1}{2}\lambda_3^2 - 2\lambda_4^2)I_0'(1/\lambda_1)}. \label{eq:pweak}
\end{align}
Motivated by Figure \ref{fig:f2}, we will now show pointwise convergence of \eqref{eq:pweak} to the classical Navier-Stokes solution, $u_0(\sigma) = 1-\sigma^2$. Clearly, $u(1) = u_0(1) = 0$, so we show pointwise convergence to the classical solution $u_0$ for $\sigma \in (0,1)$. The case $\sigma = 0$ follows a similar argument. We use the same asymptotics for $I_0(z)$ as above in addition to the following for the second derivative:
\begin{equation*}
	I_0''(z) = \frac{e^z}{(2\pi z)^{1/2}}(1+O(z^{-1})).
\end{equation*}
Using these asymptotics, we see that
\begin{align*}
	&I_0(\sigma/\lambda_1) - I_0(1/\lambda_1) = \frac{e^{\sigma/\lambda_1}\lambda_1^{1/2}}{(2\pi\sigma)^{1/2}} (1 + O(\lambda_1/\sigma)) - \frac{e^{1/\lambda_1}\lambda_1^{1/2}}{(2\pi)^{1/2}}  (1 + O(\lambda_1)), 
\end{align*}
and
\begin{align*}
	&\lambda_1^2I_0''(1/\lambda_1) - \lambda_1\left(\frac{1}{4}\lambda_2^2 + \frac{1}{2}\lambda_3^2 - 2\lambda_4^2\right)I_0'(1/\lambda_1) \\
	&= \lambda_1^2\frac{e^{1/\lambda_1}\lambda_1^{1/2}}{(2\pi)^{1/2}}(1+O(\lambda_1)) - \lambda_1\left(\frac{1}{4}\lambda_2^2 + \frac{1}{2}\lambda_3^2 - 2\lambda_4^2\right)\frac{e^{1/\lambda_1}\lambda_1^{1/2}}{(2\pi)^{1/2}}(1+O(\lambda_1)).
\end{align*}
We note that by \eqref{eq:primelength}, $\la_1 (\frac{1}{4}\la_2^2 + \frac{1}{2}\la_3^2 - 2 \la_4^2) = O(\la_1^3).$ We then conclude that
\begin{align*}
	&\frac{2\lambda_1^2(\lambda_1^2-\frac{1}{4}\lambda_2^2-\frac{1}{2}\lambda_3^2 + 2\lambda_4^2)\bigl [ I_0(\sigma/\lambda_1) - I_0(1/\lambda_1)\bigr ]}{\lambda_1^2I_0''(1/\lambda_1) - \lambda_1(\frac{1}{4}\lambda_2^2 + \frac{1}{2}\lambda_3^2 - 2\lambda_4^2)I_0'(1/\lambda_1)} \\
	&= \frac{2\lambda_1^2\left(\lambda_1^2-\frac{1}{4}\lambda_2^2 -\frac{1}{2}\lambda_3^2 + 2\lambda_4^2\right)\left(e^{-(1 - \sigma)/\lambda_1} \sigma^{-1/2} (1 + O(\lambda_1/\sigma))- (1 + O(\lambda_1))\right)}{\lambda_1^{2}(1+O(\lambda_1)) - \lambda_1\left(\frac{1}{4}\lambda_2^2 + \frac{1}{2}\lambda_3^2 - 2\lambda_4^2\right)(1+O(\lambda_1))} \\
    &\to 0
\end{align*}
as $\lambda_1 \to 0$. We have thus shown
\begin{align*}
    u(\sigma) &= 1 - \sigma^2 + \frac{2\lambda_1^2(\lambda_1^2-\frac{1}{4}\lambda_2^2-\frac{1}{2}\lambda_3^2 + 2\lambda_4^2)\bigl [ I_0(\sigma/\lambda_1) - I_0(1/\lambda_1)\bigr ]}{\lambda_1^2I_0''(1/\lambda_1) - \lambda_1(\frac{1}{4}\lambda_2^2 + \frac{1}{2}\lambda_3^2 - 2\lambda_4^2)I_0'(1/\lambda_1)} \\
    &\to 1 - \sigma^2 = u_0(\sigma) \text{ as } \lambda_1 \to 0,
\end{align*}
proving the pointwise convergence of $u(\sigma)$ to the classical solution in the weak adherence case. 

Via very similar arguments, one may also show that the dimensionless discharge rate $\Phi$ converges to the classical dimensionless rate $\Phi_0 = 1$ as $\la_1 \rar 0$ for both strong adherence and weak adherence boundary conditions; the details are left to the interested reader. 

\section{Taylor-Couette Flow}
We now consider rotational flow for a tube with an outer radius $R$ rotating with angular velocity $\Omega$. The velocity and pressure are taken to be of the form
\begin{align}
	\bs v = v(r)\bs e_\theta, \quad p = p(r). 
\end{align}
Then 
\begin{align}
	\nabla \bs v &= v' \bs e_\theta \tens \bs e_r - \frac{1}{r} v \bs e_r \tens \bs e_\theta, \\
	\nabla \nabla \bs v &= v'' \bs e_\theta \tens \bs e_r \tens \bs e_r - \Bigl ( \frac{v'}{r} - \frac{v}{r^2} \Bigr ) \bs e_r \tens \bs e_r \tens \bs e_\theta \\
	&\quad  - \Bigl ( \frac{v'}{r} - \frac{v}{r^2} \Bigr ) \bs e_r \tens \bs e_\theta \tens \bs e_r + 
	 \Bigl ( \frac{v'}{r} - \frac{v}{r^2} \Bigr ) \bs e_\theta \tens \bs e_\theta \tens \bs e_\theta.
\end{align}

\subsection{General solution}
The $\theta$-component of \eqref{eq:field} for $\bs v = v(r)\bs e_\theta$ reads
\begin{equation}
  0 = \cl L
  \left(1 - \ell_1^2\cl L\right)v
\end{equation}
where $\cl L = \frac{d^2}{dr^2} + \frac{1}{r}\frac{d}{dr}-\frac{1}{r^2}$. The above equation has the following general solution:
\begin{equation}
  v(r) = \Omega r + c_1 r + c_2 \frac{1}{r} + d_1 I_1(r/\ell_1) + d_2 K_1(r/\ell_1)
\end{equation}
with constants $c_1, c_2, d_1$ and $d_2$ determined by the boundary conditions. We require $\bs v$ to be regular at the origin and thus, $c_2 = d_2 = 0$. We now determine the solutions satisfying either strong adherence or weak adherence boundary conditions.

\subsection{Strong adherence boundary conditions}
The boundary condition $\p_{\bs n} \bs v = \bs{0}$ is equivalent to the requirement that $v'(R) = 0$. Thus, the strong adherence boundary conditions are as follows:
\begin{align*}
	\Omega R &= \Omega R + c_1R + d_1 I_1(R/\ell_1) \\
	0 &= \Omega + c_1 + d_1 I_1(R/\ell_1).
\end{align*}
It is readily checked that the solutions are
\[
c_1 = -\frac{\Om R I_1(R/\ell_1)}{I_1(R/\ell_1) - \frac{R}{\ell_1}I_1'(R/\ell_1)},\quad d_1 =\frac{\Om R}{I_1(R/\ell_1) - \frac{R}{\ell_1}I_1'(R/\ell_1)}
\]
so the solution can be written as
\[
v(r) = \Om r + \frac{\Om}{I_1(R/\ell_1) - \frac{R}{\ell_1}I_1'(R/\ell_1)}(RI_1(r/\ell_1) - rI_1(R/\ell_1)).
\]
\begin{figure}[t]
	\centering
	\includegraphics[trim = {0 6cm 0 6cm},clip,scale = 0.55]{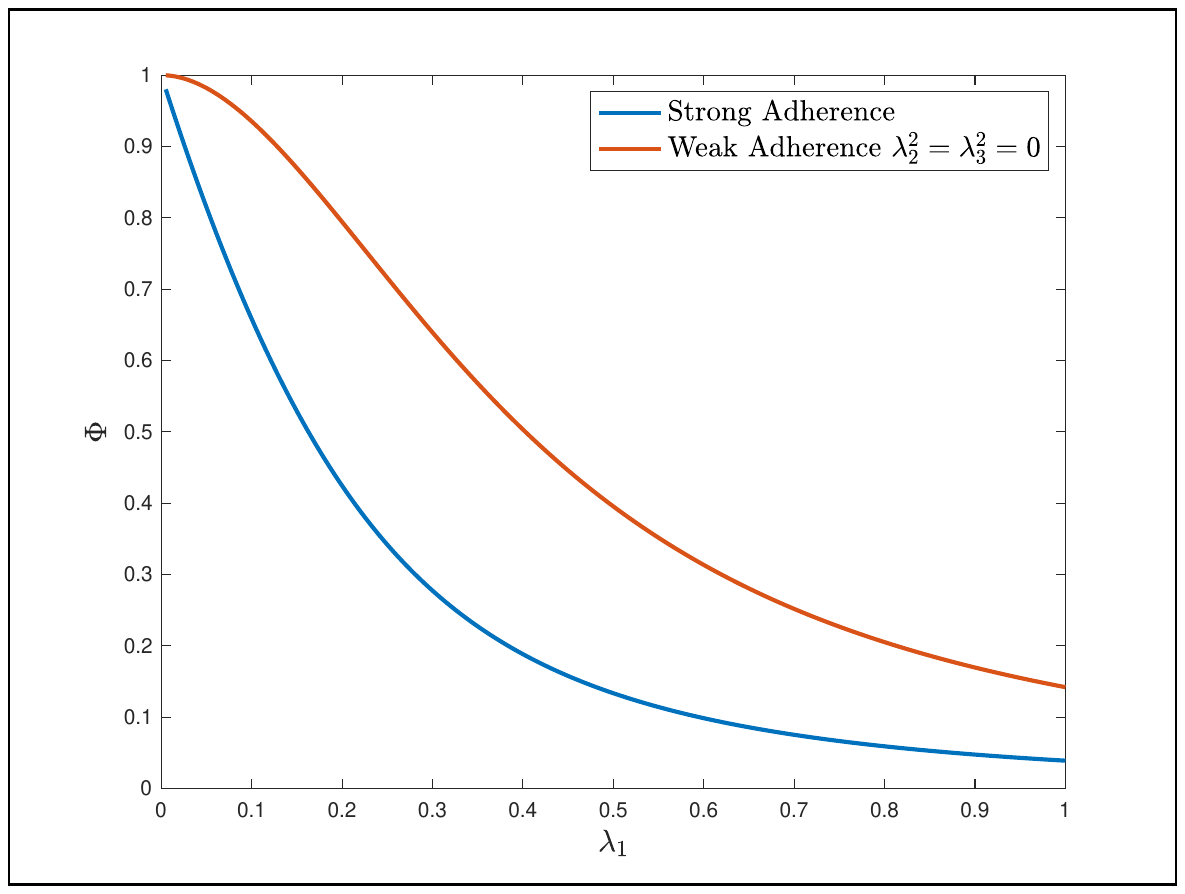}
	\caption{Graph of $\Phi$ with strong adherence boundary conditions and weak adherence boundary conditions, where $\lambda_2 = \lambda_3 = 0$.}
	\label{fig:f4}
\end{figure}
The dimensionless velocity $u(\sigma) = (\Omega R)^{-1} v(R\sigma)$ as a function of $\sigma = r/R$ is given by 
\begin{align}
u(\sigma) = \sigma + \frac{I_1(\sigma/\lambda_1)-\sigma I_1(1/\lambda_1)}{I_1(1/\lambda_1) - \lambda_1^{-1}I_1'(1/\lambda_1)}, \label{eq:dimensionlessstrongTC}
\end{align}
where $\lambda_1 = \ell_1/R$.

\subsection{Weak adherence boundary conditions}

For the weak adherence case, we assume that $\bbG [\bs n \tens \bs n] = \bs 0$. Since we assume that $p'(R) = 0$, weak adherence requires that
\begin{align}
	(\eta_1 - \eta_2 - 4\eta_3) v''(R) - (3 \eta_2 + 4 \eta_3)\Bigl ( \frac{1}{R}v'(R) - \frac{1}{R^2}v(R) \Bigr ) = 0. \label{eq:weaktc} 
\end{align}
We observe that the classical solution $v(r) = \Omega r$ satisfies $v(R) = \Omega R$ and \eqref{eq:weaktc} {$I_1(r)$ does not}, and thus, the velocity profile satisfying weak adherence conditions is $$v(r) = \Omega r.$$ The nondimensionalized solution is then given by 
\begin{equation}
	u(\sigma) = \sigma.
\end{equation} 

\begin{figure}[t]
	\centering
	\includegraphics[trim = {0 6cm 0 6cm},clip,scale = 0.55]{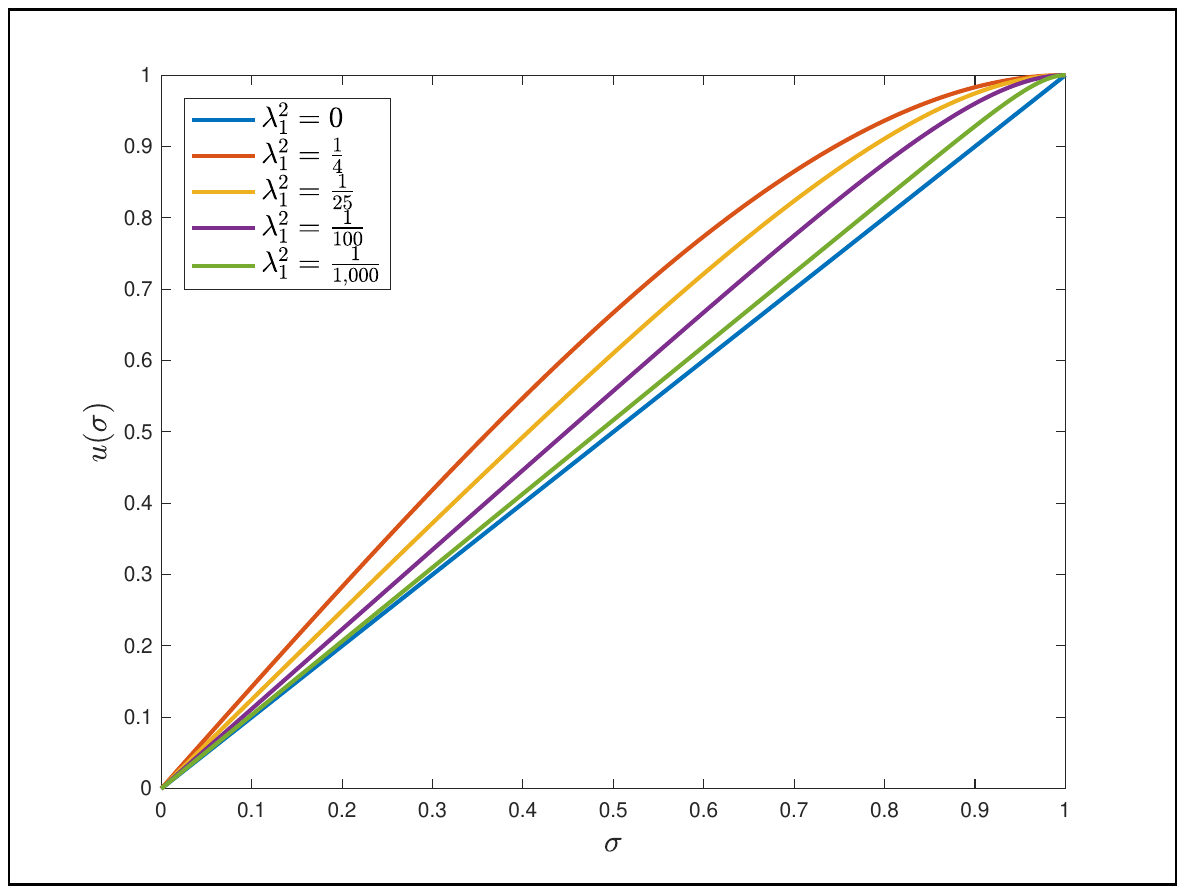}
	\caption{Graph of $u(\sigma)$ for Couette flow and strong adherence boundary conditions. The value $\lambda_1 = 0$ corresponds to the classical solution.}
	\label{fig:f5}
\end{figure}
\subsection{Convergence to the solution for the classical Navier-Stokes model}
Motivated by Figure \ref{fig:f5}, we now prove pointwise convergence of the dimensionless velocity field $u(\sigma)$ given by \eqref{eq:dimensionlessstrongTC} to the dimensionless velocity field corresponding to the classical Navier-Stokes model, $u_0(\sigma) = \sigma$. 

Clearly $u(0) = u_0(0)$ and $u(1) = u_0(1)$, so, we only consider $\sigma \in (0,1)$. We use the asymptotics for $I_1(z)$ \cite{abramowitz1964handbook}: 
\begin{align}
	I_1(z) = \frac{e^z}{(2\pi z)^{1/2}}\Bigl (1 + O(z^{-1}) \Bigr ), \quad 
	I_1'(z) = \frac{e^z}{(2\pi z)^{1/2}} \Bigl ( 1 + O(z^{-1}) \Bigr ),
\end{align}
as $z \rar \infty$. Then for $\sigma \in (0,1)$, 
\begin{align}
	I_1(\sigma/\lambda_1) - \sigma I_1'(1/\lambda_1) &=
	(\lambda_1/2\pi \sigma)^{1/2} e^{\sigma/\lambda_1} ( 1 + O(\lambda_1/\sigma) ) \\&\quad - \sigma ( \lambda_1/2\pi)^{1/2} e^{1/\lambda_1}(1 + O(\lambda_1)), \\
	I_1(1/\lambda_1) - \lambda_1^{-1} I_1'(1/\lambda_1) &= {(\lambda_1/2\pi)^{1/2}e^{1/\lambda_1}(1 + O(\lambda_1))}\\
	&-(1/2\pi\lambda_1)^{1/2} e^{1/\lambda_1}(1 + O(\lambda_1)).
\end{align}
We conclude that for each $\sigma \in (0,1]$, 
\begin{gather}
\frac{I_1(\sigma/\lambda_1)-\sigma I_1(1/\lambda_1)}{I_1(1/\lambda_1) - \lambda_1^{-1}I_1'(1/\lambda_1)} = -\frac{\lambda_1 \sigma^{-1/2}e^{-(1-\sigma)/\lambda_1}(1 + O(\lambda_1/\sigma)) - \sigma \lambda_1(1 + O(\lambda_1))}{{(\lambda_1 - 1)}(1+O(\lambda_1))} \\
\rar 0, \quad \mbox{as } \lambda_1 \rar 0. 
\end{gather}
Thus, for all $\sigma \in [0,1]$, 
\begin{align}
\sigma + \frac{I_1(\sigma/\lambda_1)-\sigma I_1(1/\lambda_1)}{I_1(1/\lambda_1) - \lambda_1^{-1}I_1'(1/\lambda_1)} \rar \sigma = u_0(\sigma), 
\end{align}
as $\lambda_1 \rar 0$, proving pointwise convergence to the solution for the classical Navier-Stokes model. Since the classical solution and the weak adherence solution are the same, there is no convergence to show for weak adherence boundary conditions.

\subsection{The pressure}
The $r$-component for Taylor-Couette flow determines the pressure, which we assume is regular and satisfies $p'(R) = \frac{\p p}{\p \bs n}(R) = 0$. The pressure $p(r)$ satisfies the following differential equation:
\begin{equation} \label{eqn:TCPressure}
  \Bigl ( p'' + \frac{1}{r}p' - \frac{1}{\ell_1^2}p\Bigr )' = \rho \frac{\ell_0^2}{\ell_1^2}\Bigl (v \Bigl (\frac{v(r)}{r}\Bigr )' + \Bigl (\frac{v(r)}{r}\Bigr )^2\Bigr )' - \frac{\rho}{\ell_1^2} \frac{(v(r))^2}{r},
\end{equation}
where $v(r)$ is the velocity. Carrying out the differentiation, we see that $p'$ satisfies 
\begin{align}
	(p')'' + \frac{1}{r}(p')' - \Bigl (\frac{1}{\ell_1^2} + \frac{1}{r^2} \Bigr ) p' = \rho \frac{\ell_0^2}{\ell_1^2}\Bigl (v \Bigl (\frac{v(r)}{r}\Bigr )' + \Bigl (\frac{v(r)}{r}\Bigr )^2\Bigr )' - \frac{\rho}{\ell_1^2} \frac{(v(r))^2}{r} 
\end{align} 

We will label the right-hand side as $f(r)$. Using the change of variables $h(r) = r^{1/2}p'(r)$, we have that the ODE that $h$ satisfies is
\begin{equation} \label{eqn:HPressureODE}
 \begin{cases} 
  h''(r) + \Bigl (\frac{1}{4r^2} - \Bigl (\frac{1}{r^2} + \frac{1}{\ell_1^2}\Bigr )\Bigr ) = r^{1/2}f(r), &  r \in (0,R) \\
  h(0) = h(R) = 0.
 \end{cases}
\end{equation}
\begin{figure}[t]
	\centering 
	\includegraphics[trim = {0 6cm 0 6cm},clip,scale = 0.55]{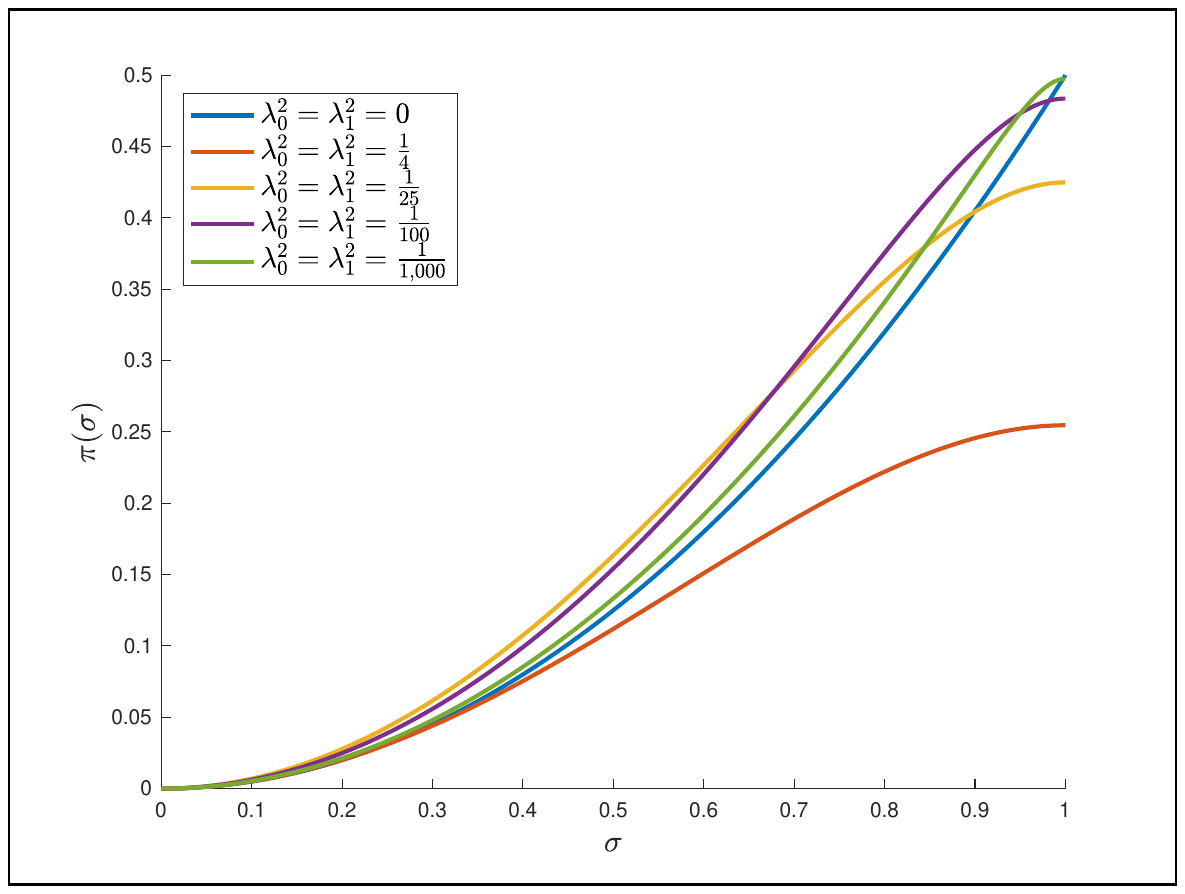}
	\caption{Graph of the dimensionless pressure $\pi(\sigma):=\int_0^\sigma \pi(s)ds$ for Couette flow and strong adherence boundary conditions. The curve corresponding to $\lambda_0 = \lambda_1 = 0$ is the classical solution. The Numerical solutions for $\pi$ were obtained using \textsc{bvp4c} \textsc{MATLAB} solver.} 
	\label{fig:f6}
\end{figure}

The solution to this ODE is given as the following via variation of parameters:
\begin{align*}
  h(r) &= r^{1/2}I_1(r/\ell_1 ) \int_r^{R} \Bigl ( \frac{K_1(R/\ell_1)}{I_1(R/\ell_1)}I_1(s/\ell_1) - K_1(s/\ell_1) \Bigr )sf(s)\,{d}s \\
  &+ r^{1/2} \Bigl (\frac{K_1(R/\ell_1)}{I_1(R/\ell_1)}I_1(r/\ell_1) - K_1(r/\ell_1)\Bigr )\int_0^r I_1(s/\ell_1)sf(s)\,{d}s.
\end{align*}

Thus,
\begin{align*}
  p'(r) &=  I_1(r/\ell_1 ) \int_r^{R} \Bigl ( \frac{K_1(R/\ell_1)}{I_1(R/\ell_1)}I_1(s/\ell_1) - K_1(s/\ell_1) \Bigr )sf(s)\,{d}s \\
  &+ \Bigl (\frac{K_1(R/\ell_1)}{I_1(R/\ell_1)}I_1(r/\ell_1) - K_1(r/\ell_1)\Bigr )\int_0^r I_1(s/\ell_1)sf(s)\,{d}s   
\end{align*}
With $\phi(\sigma) = \frac{R}{\rho \Omega^2} f(R \sigma)$, the dimensionless pressure $\pi(\sigma) = (\rho \Omega^2 R^2)^{-1} p(R \sigma)$ with $\sigma := r/R$ then satisfies 
\begin{align}
	\pi'(\sigma) &= I_1(\sigma/\lambda_1) \int_\sigma^1 \Bigl (  \frac{K_1(1/\lambda_1)}{I_1(1/\lambda_1)}I_1(s/\lambda_1) - K_1(s/\lambda_1) \Bigr ) s \phi(s) ds \\
	&+ \Bigl (\frac{K_1(1/\la_1)}{I_1(1/\la_1)}I_1(\sigma/\la_1) - K_1(\sigma/\la_1) \Bigr )\int_0^\sigma I_1(s/\la_1)s\phi(s)\,{d}s. 
\end{align}


As shown in Figure \ref{fig:f6}, the dimensionless pressure field approaches the classical solution. A comprehensive theory establishing convergence of both the nondimensional velocity and pressure fields to their classical limits will be pursued elsewhere.

\bibliographystyle{plain}
\bibliography{researchbibmech, researchbibmech_NSF}

\bigskip

\centerline{\scshape C. Balitactac}
\smallskip
{\footnotesize
	\centerline{Department of Mathematics, University of North Carolina}
	
	\centerline{Chapel Hill, NC 27599, USA}
	
	\centerline{\email{corbindb@unc.edu}}
}

\bigskip

\centerline{\scshape C. Rodriguez}
\smallskip
{\footnotesize
	\centerline{Department of Mathematics, University of North Carolina}
	
	\centerline{Chapel Hill, NC 27599, USA}
	
	\centerline{\email{crodrig@email.unc.edu}}
}

\end{document}